\documentclass[a4paper,11pt,reqno]{amsart}
\usepackage{amsthm}
\usepackage{amsmath}
\usepackage{enumerate}
\usepackage{amsfonts}
\usepackage{amssymb}
\usepackage{fullpage}
\usepackage{amsmath,amscd}
\usepackage{stmaryrd}
\usepackage{graphicx}
\usepackage{array}
\usepackage{amsmath,amssymb,amsfonts,dsfont}
\usepackage[utf8]{inputenc} 
\usepackage[english]{babel}
\usepackage[T1]{fontenc} 
\usepackage{graphicx}
\usepackage{float}
\usepackage{pdfpages}
\usepackage[a4paper, margin = 3cm, bottom = 3cm]{geometry}
\usepackage{ifpdf}
\usepackage{marginnote}
\usepackage{hyperref}	
\usepackage{pict2e}
\usepackage{epstopdf}
\usepackage{epsfig,graphicx,amsmath,amsfonts}
\usepackage{xcolor,import}
\hypersetup{pdfborder=0 0 0, 
	    colorlinks=true,
	    citecolor=black,
	    linkcolor=blue,
	    urlcolor=red,
	    pdfauthor={Mathieu Kohli}
	   }
\newtheorem{thm}{Theorem}[section]

\newtheorem{algo}[thm]{Algorithm}
\theoremstyle{definition}
\newtheorem{defn}[thm]{Definition}
\theoremstyle{remark}

\theoremstyle{definition}

\theoremstyle{definition}

\theoremstyle{definition}

\numberwithin{equation}{section} 
\title{A mathematical modelization of the prosperity of an Apis Mellifera beehive}
\author{Mathieu Kohli} 
\email{Mathieu.Kohli@protonmail.ch}
\begin{document}
%
%
%
\maketitle
\setcounter{tocdepth}{1}
\tableofcontents
%
%
\section*{Introduction}

Aesop's famous fable \emph{The Ant and the Grasshopper} emphasizes the choice that living beings have to make between survival (the Ant)  and reproduction (the Grasshopper). Apis Mellifera colonies face the same dilemma : must they rather collect nectar in order to produce honey that will enable them to endure the cold of winter or should they choose to bring pollen back to the hive so as to raise young bees ? Both, of course, but in which proportion ? And where should they fly to find the best resources ?
\newline
We adress these question through a mathematical perspective.
\newline
We only consider healthy colonies and we assume that they optimize the use of resources, which we already know they do in a specific case through the honeycomb theorem \cite{HoneycombTheorem}.

\section{Energy, matter and bees}

To "build" one bee, the colony must feed the larva a certain weight of pollen (matter) and a certain weight of honey (energy). Once it leaves from the nest and until it dies, the bee does not undergo important morphological transformation and we therefore assume as a first approximation that the quantity of pollen and honey it "contains" is constant. Clearly, the bees continue to eat after they leave the nest but this is only done to compensate the energy that is naturally lost through their life process and to fight the material decay of the body, which means that bees eat to maintain a constant body. In other words, if we only focus on energy and matter, a bee can be considered as a box in which a constant quantity of pollen and of honey is hidden.

\subsection{Energy}
Inside the hive, there are three forms of energy : honey, bees and wax.
\begin{equation*}
E = \mu M + \alpha N + \gamma C,
\end{equation*}
where $E$ represents energy, $M$ and $C$ are the weight of honey and honeycombs respectively contained in the hive, $N$ stands for the number of bees in the colony (which incudes the foragers that are outside the hive).
\newline
Since we need to record precise information, $N$ may be seen as an element of $\mathbb{R}^\mathcal{L}\times \mathbb{R}$ where $\mathcal{L}$ denotes the maximal lifespan (counted in days) of a female bee. The component $N_i$ contains the number of bees whose age is $i$ and $N_{\mathcal{L}+1}$ equals the number of males inside the hive.
\newline
For notations to be coherent, we consider $\alpha\in\left(\mathbb{R}^{\mathcal{L}}\times\mathbb{R}\right)^*$. The assumption that each bee "contains" a constant energy translates into a covector $\alpha$ whose $\mathcal{L}$ first components contain a fixed constant.
\newline
We can now study the temporal evolution of the energy in the hive. First we must define what we mean by "the hive". We see "the hive" as a blackbox encompassing the pollen and honey reserves, as well as the wax honeycombs and all the living bees that belong to the colony (even the foragers that are flying miles away from the other bees). We have :
\begin{equation}\label{bilanenergie}
\boxed{\frac{\text{d} E}{\text{d}t}=\mu\frac{\delta M_{foraging}}{\text{d}t} -\alpha \frac{\delta N_{death}}{\text{d}t}  - \pi N -h \left(T, N, N_{males}, N_{nurses}, N_{larvae}  \right),}
\end{equation}
In the previous energetic "balance sheet" :
\newline $\bullet$ $\frac{\delta M_{foraging}}{\text{d}t}$ stands for the quantity of honey produced with the nectar entering the hive to which we substract the honey that has been used in order to transform this nectar into honey.
\newline $\bullet$  $h$ is the cost of heating the hive.
\newline $\bullet$  $\pi N$ represents the "households expenditures" of the hive.
\subsection{Matter} A equation similar to \ref{bilanenergie} can be written for pollen, except that there is no term that corresponds to $h$ in the pollinic equation.
\subsection{Size of the colony}We defined $N_j$ as the number of bees which are $j$ days old.
We also know that the task of a bee in the colony depends on its age. If bees perform their $i$-th task between age $D_i$ and age $D_{i+1}-1$, we can compute the number of bees performing a specific task as
\begin{equation*}
\tilde{N}_i=\sum_{d=D_i}^{D_{i+1}-1}N_d.
\end{equation*}

We consider the survival rate
\begin{equation*}
\begin{matrix}
s&:&\left[\left|0,\mathcal{L} \right|\right]&\longrightarrow &\left[ 0,1 \right],
\end{matrix}
\end{equation*}
in which we read the proportion of bees that survive at least $i\in \left[\left|0,\mathcal{L} \right|\right]$ days. The survival rate depends on the quality of the pollen with which the bees were fed when they were larvae. \cite{PollenNutrition}
\newline
At time $t$,
\begin{equation*}
\tilde{N}_i\left( t\right)=\sum_{d=D_i}^{D_{i+1}-1}N_0\left(t-d\right)s(d).
\end{equation*}
Moreover, we express the mortality term in $\ref{bilanenergie}$  as
\begin{equation*}
\boxed{\frac{\delta N_{death}}{\text{d}t}=-\sum_{d=0}^{\mathcal{L}}N_d \left.\frac{\text{d}s}{\text{d}t}\right|_{t=d}.}
\end{equation*}

\section{Heating the hive}
There are two main reasons for bees to heat the hive :
\newline $\bullet$ to maintain an ideal temperature $T_{brood}\simeq 35.5^\circ C$ for brood cells \cite{Thermoregulation1} which is part of the nurses' job. In order to do this, $\frac{N_{nurses}}{N_{larvae}}$ needs to be high enough for the nurses to take care of the whole brood. If this is not the case, some of the brood cells will be abandonned.
\newline $\bullet$ to have the temperature at the center of the hive above a minimal value $T_{center}$ even in winter when there are no brood cells \cite{Thermoregulation1}.
\newline
We study thermoregulation through the heat equation and also refer to works such as \cite{Thermoregulation1,Thermoregulation2}.
\newline
To reach the right temperature, different strategies can be implemented :
\newline $\bullet$ If the hive is active, no compact mass of bees can be formed and each bee has to heat its own surroundings so that the power spent by the hive is proportional to the number $N_h$ of bees that are taking part in the heating process and to the difference between the targeted temperature $T_{target}$ and the temperature outside the hive $T_{out}$ :
\begin{equation*}
h\left( N_h, T_{out} \right) = \theta N_h \left| T_{target}-T_{out} \right|.
\end{equation*}
\newline $\bullet$ When the hive is inactive, which happens during winter for instance, the bees can choose a more efficient mechanism to heat the hive. They form a compact mass where each bee keeps as close as possible to its immediate neighbours. Instead of heating only itself, each bee contributes to heating the bee cluster. We apply the heat equation to the cluster as a whole to find out how much power is spent to stabilize the temperature at $T_{target}$ at the center of the cluster. For the sake of simplicity, we assume that the bee cluster is a $R$-radius sphere. The heat equation reads :
\begin{equation}
\frac{1}{r^2}\frac{\partial}{\partial r}\left( r^2 \frac{\partial T}{\partial r} \right) +Q(r)=\overbrace{\frac{\partial T}{\partial t} = 0}^{\text{stationarity}}.\label{EquationChaleur}
\end{equation}
Let us notice that $T$ is an even analytic function such that $T(0)=T_{target}$ and $T(R)=T_{out}$ so that we can write
\begin{equation*}
T(r)=T_{target}-\alpha r^{2\nu}+ o\left(r^{2\nu}\right),
\end{equation*}
and it follows from the heat equation and the assumption that $T$ is analytic that

\begin{align*}
Q(r)&= 2\nu\left( 2\nu+1 \right)\alpha  r^{2\nu-2}+ o\left(r^{2\nu-2}\right),\\
&\simeq  2\nu\left( 2\nu+1 \right) \frac{T(r)-T_{target}}{r^2}, \\
&\simeq  2\nu\left( 2\nu+1 \right) \left(T(r)-T_{target}\right) \left(\frac{T'(r)}{2\nu\left(T(r)-T_{target}\right)}\right)^2
\end{align*}
and
\begin{equation}
Q(r)\simeq  \frac{2\nu+1}{2\nu}  \frac{\left\Vert \overrightarrow{\text{grad}}T \right\Vert^2}{T(r)-T_{target}} \text{ for $r\simeq 0$.}\label{ChaleurRegal0}
\end{equation}
But a bee inside the bee cluster must "compute" the heat power it emits with information it percieves from its surroundings only. In particular, the bee does not know its exact position $r$ inside the cluster and the formula it uses to compute the heat it emits must not depend on $r$. Since formula \ref{ChaleurRegal0} is computed by using local information only, it is a valid formula for bees to compute the heat they emit and since bees must use the same formula at any position $x$ inside the cluster, we can write :
\begin{equation}
\boxed{Q(x) =\frac{2\nu+1}{2\nu}  \frac{\left\Vert \overrightarrow{\text{grad}}T \right\Vert^2}{T(x)-T_{target}}. }
\end{equation}
We notice that the previous formula makes sense even when there are no spherical symetries, but we keep these symetries in order to make computations easier.
\newline
Now we know the power each bee emits, we can find out what the temperature is inside the cluster and how much power is used by the entire bee cluster. The following temperature satisfies \ref{EquationChaleur}, is even and such that $T(0)=T_{target}$, $T(R)=T_{out}$ :
\begin{equation}\label{TemperatureCluster}
\boxed{T(r)=T_{target}-\left( T_{target}-T_{out} \right)\left( \frac{r}{R} \right)^{2\nu}.}
\end{equation}
Now if we have $N$ bees in the cluster, $R$ is proportional to $N^{1/3}$ and the power the bee colony spends to heat itself is
\begin{align*}
h(N,T_{out})&= \int_{0}^{R(N)}\int_{\theta}\int_{\phi}Q(r)r^2\text{d}r\text{d}\theta\text{d}\phi, \\
&=\int_{0}^{R(N)}Q(r)r^2\text{d}r\int_{\theta}\text{d}\theta\int_{\phi}\text{d}\phi,\\
&=\int_{0}^{R(N)}2\nu\left( 2\nu+1 \right)\left( T_{target}-T_{out} \right)\left( \frac{r}{R(N)} \right)^{2\nu}\text{d}r\int_{\theta}\text{d}\theta\int_{\phi}\text{d}\phi, \\
&=\int_{0}^{1}2\nu\left( 2\nu+1 \right)\left( T_{target}-T_{out} \right)x^{2\nu}R(N)\text{d}x\int_{\theta}\text{d}\theta\int_{\phi}\text{d}\phi.
\end{align*}
The previous expression is proportional to $R(N)\left( T_{target}-T_{out} \right)$ so we can conclude
that
\begin{equation}
\boxed{h(N,T_{out})=\mathcal{K} N^{1/3}\left( T_{target}-T_{out} \right).}
\end{equation}
\section{Available floral resources}
In order to bring nectar back to their hive, foragers have to reach a floral resource, fill themselves with as much nectar as possible and come back to the hive. If we consider a benchmark flower spieces whose nectar quality is optimal, we denote by $q_0$ the quantum of energy contained in the nectar that one bee can transport to which we substract the energy that has to be spent by the hive to transform this quantity of nectar into honey (which includes for example the cost of evacuating humidity from the hive).
\newline
For each trip a forager collecting nectar makes to this benchmark spieces, the hive benefits from an energetical surplus whose value is
\begin{equation*}
q_0-\sigma d - \text{ foraging cost,}
\end{equation*}
where $d$ is the distance between the hive and the floral resource.
\newline
There exists a distance $d_{max}$ $(\text{more than } 10km)$ beyond which no forager bee will venture to collect nectar \cite{MaxDist,MaxDistance}. This distance $d_{max}$ corresponds to the limit beyond which the energetical benefit of foraging is zero so that when $q_0$ is associated to the benchmark flower with optimal intrinsic quality,
\begin{equation*}
q_0-\sigma d_{max}=0 \Longrightarrow \sigma =\frac{q_0}{d_{max}}.
\end{equation*}
The energetical benefit the hive gets from one bee foraging at distance $d$ on a flower $f$ (not necessarily the benchmark flower) is
\begin{equation}\label{CoutDistance}
q_f - q_0 \frac{d}{d_{max}}-\text{ foraging cost.}
\end{equation}
Now we must compute the foraging cost. There are two different parts in the foraging process : the bee is either on a flower collecting nectar or hopping from one flower to another. Since foraging is a symbiotic activity with the flower, it is in the best interest of the flower to minimize the inconvenience for the bee to visit it. However, it is in the best interest of the flowers to maximize the wandering of the bees between flowers to spread their pollen, and since they compensate the bees with nectar, bees "agree" to this costly process.
\newline
We ignore the cost of visiting a flower which both flowers and bees "try" to minimize and focus on the cost of hopping from one flower to another.
\newline
During one trip to a floral field, each bee has to visit a certain number $m_f$ of flowers ($f$ is the flower spieces) in order to be filled with nectar. The average distance between each flower is
\begin{equation*}
\overline{d}=\left\lbrace\begin{matrix}
\frac{k_2}{\sqrt{\rho_f}} &\text{ if the flowers grow in a field,}\\
\frac{k_3}{\sqrt[3]{\rho_f}} &\text{ if the flowers are on a bush or a tree,}
\end{matrix}\right.
\end{equation*}
where $\rho_f$ corresponds to the floral density.
\newline
Moreover, we can measure $v_{hop}$, the average speed of the bees hopping between flowers that is different from $V$, the speed at which bees travel long distances. Since we already know thanks to \ref{CoutDistance} how much it costs to travel at speed $V$ on a distance $d$, we can deduce the cost of hopping between flowers is equal to
\begin{equation*}
q_0\frac{1}{\sqrt[n]{\rho_f}}\frac{k_n m_f v_{hop}}{d_{max}V}.
\end{equation*}
Indeed the power dissipated during flight is proportional to the speed of the bee according to \cite{Aerodynamic,Aerodynamic2}. The exponent of $\frac{v_{hop}}{V}$ can be adapted depending on the type of flight. The specific flying pattern of bees has been studied in \cite{FlightDescription}.
\begin{defn} 

The critical floral density for flower $f$ naturally arises as
\begin{equation*}
\rho_{crit,f} = \left( \frac{k_n m_f v_{hop}}{d_{max}V} \right)^n,
\end{equation*}
where $n=2$ if the floral resource is a surface such as a field and $n=3$ if the floral resource is a tree or a bush.
\end{defn}
The critical density enables us to actually write the benefit the hive gets each time one bee forages at distance $d$ from the hive on flower $f$ :
\begin{equation*}
q_0\underbrace{\left( \frac{q_f}{q_0} -\frac{d}{d_{max}}- \sqrt[n]{\frac{\rho_{crit,f}}{\rho_f}} \right)}_{\mathcal{Q}_f}.
\end{equation*}
\begin{defn} \label{DefQualityNectar} The quality field of a floral resource $f$ is a scalar field which maps the value of the best yield the hive can obtain from this floral resource to the geographical point $x$ where the hive sits
\begin{equation*}
\begin{matrix}
\mathcal{Q}_f &: & \mathbb{R}^2 & \longrightarrow & \left[0,1\right] \\
&& x & \longmapsto & \max_{y\in\mathbb{R}^2} \left( \frac{q_f}{q_0} -\frac{d(x,y)}{d_{max}}- \sqrt[n]{\frac{\rho_{crit,f}}{\rho_f(y)}} \right),
\end{matrix}
\end{equation*}
where $n$ is the mathematical dimension of the floral resource.
\end{defn}
In a floral vacuum, the following local equation is satisfied
\begin{equation}
\boxed{\left\Vert \overrightarrow{\text{grad}} \mathcal{Q}_f\right\Vert = \frac{1}{d_{max}}.}
\end{equation}
This local equation has a very simple geometric interpretation : if on a map, we know a line of constant quality, let us say the edge of a field whose quality is 0.7, if we want to find the line of constant quality 0.6, we "roll a wheel" with diameter $\left(0.7-0.6\right)d_{max}$ on the edge of this field and the furthest points from the field that are reached by the wheel form the line of constant quality 0.6.

\begin{center}
\includegraphics[scale=.5]{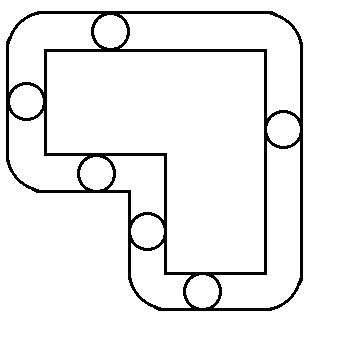}
\end{center}
\section{Appointing bees to a floral resource}
Von Frisch discovered that bees communicate the position of a floral resource using a dancing pattern \cite{Dance}. If a bee discovers a high quality resource it dances to recruit the right number of foragers to fetch the nectar that is beeing produced.
\newline
The number of bees $\phi_f$ that forage on a flower $f$ is optimal when
\begin{equation*}
\frac{\phi_f}{t_{hive}+t_{flight}+t_{foraging}} =S\rho_{f}\frac{\lambda_{f}}{q_f},
\end{equation*}
where $S$ is the surface of the floral resource, $\lambda_f$ is the surfacic power offered to bees by flower $f$ and where
\begin{equation*}
t_{hive}+t_{flight}+t_{foraging}=t_{hive}+2\frac{d}{V}+m_f\left(\beta_f+\frac{k_n}{v_{hop}\sqrt[n]{\rho_f}} \right)
\end{equation*}
is the temporal period of the cyclic movement back and forth between the hive and the flower. Experimental data has been gathered and can be found for instance in \cite{NectarFlower,RFID,RFID2}.
\newline
To summarize, we write the number of bees that visit a floral resource $f$ as
\begin{equation}\label{NumberBeesresource}
\boxed{\phi_f = S\rho_{f}\frac{\lambda_{f}}{q_f}\times \left( t_{hive}+2\frac{d}{V}+m_f\left(\beta_f+\frac{k_n}{v_{hop}\sqrt[n]{\rho_f}} \right) \right).}
\end{equation}
Every term in the previous formula has already been explicitely defined except $t_{hive}$ which is the average time a forager spends in the hive after each trip to the floral resource and $\beta_f$ which is the average time a forager spends on one individual flower $f$.
\newline
If the optimal number of bees is sent to a floral resource $f$ by the beehive, the power that the hive gains and that corresponds to an additive term in \ref{bilanenergie} is
\begin{equation}\label{Quantityresource}
\boxed{\mu\frac{\delta M_{forag,f}}{\text{d}t}=S\rho_f\lambda_f\frac{q_0}{q_f}\mathcal{Q}_f.}
\end{equation}
A central notion for bees arises from combining the above formulas. The hive's efficiency on a floral resource corresponds to the quantity of resources one bee brings back to the hive every time unit, and that is computed by dividing \ref{Quantityresource} by \ref{NumberBeesresource}.
\begin{defn} \label{DefEfficient} The efficiency of the hive on a floral resource $f$ is
\begin{align*}
\eta_{f}=\frac{q_0 \mathcal{Q}_f}{t_{hive}+2\frac{d}{V}+m_f\left(\beta_f+\frac{k_n}{v_{hop}\sqrt[n]{\rho_f}} \right)}.
\end{align*}
\end{defn}
\section{Pollen/Honey exchange rate}
Until now we have only focused on flowers that offer nectar to bees. What about pollen ? How does the hive choose to fetch either nectar or pollen ? We must come back to the dilemma of \emph{the Ant and the Grasshopper}. The main goal of the colony is to survive winter. In order to achieve this goal, the beehive must enter the cold season with an optimal ratio between matter and energy that can be "decomposed" as the ratio of energy to matter that is used to form each bee to which we add the ratio of the energy that each bee will need to spend during winter to the pollinic matter that was used to build the bee.
\begin{equation*}
\mathcal{R}_{target}= \left. \frac{\text{Honey}}{\text{Pollen}}\right\vert_{\text{bee}} +\frac{\int_{winter}\mathcal{K} N^{1/3} \left( T_{target}-T(t) \right)  \text{d}t }{N \times \left.\text{ Pollen}\right\vert_{\text{one bee}}}.
\end{equation*}
Notice that the target ratio depends on the number of bees $N$ that will constitute the winter colony and that this ratio decreases with $N$. To ensure prosperity the colony should consider the worst case scenario where $N$ is small as it fixes a target ratio $\mathcal{R}_{target}$. When winter comes closer, $N$ can be predicted more acurately and $\mathcal{R}_{target}$ can be defined more precisely.
\newline
At any moment the hive happens to be in one of three situations :
\newline $\bullet$ $\mathcal{R}_{hive}<\mathcal{R}_{target}$, which happens for instance after winter when the honey has been depleted,
\newline $\bullet$ $\mathcal{R}_{hive}>\mathcal{R}_{target}$,
\newline $\bullet$ $\mathcal{R}_{hive}\simeq\mathcal{R}_{target}$, which should be the case before winter and also happens between the two above situations as a transition phase. Such a transition implies that the bees need to be reallocated to new floral resources which takes time as bees need to explore new places and new foragers have to replace old foragers. Moreover, the size of the nest is limited by the queen that can only lay eggs at a certain speed and by the number of nurses that tend to the nest at some point in time.
The limitation to the growth of the nest by the queen is linear and the limitation that comes from the number of nurses is exponentially lifted as eggs soon become nurses. Since the size of the nest is limited by these factors, as the hive transitions from $\mathcal{R}_{hive}<\mathcal{R}_{target}$ to $\mathcal{R}_{hive}>\mathcal{R}_{target}$, more pollen will be needed, but this will be a progressive  continuous change and not a brutal change which will lead $\mathcal{R}_{hive}$ to overshoot. This overshot is important for beekeepers as they use it to collect honey.
\newline
We now study the two cases $\mathcal{R}_{hive}<\mathcal{R}_{target}$ and $\mathcal{R}_{hive}>\mathcal{R}_{target}$.
\newline
Our goal is to define a quantity $\tau$ called the exchange rate, that shall represent the "price" of pollen paid in honey.
\newline
We already defined the quality of a flower that produces nectar in Definition \ref{DefQualityNectar}. For pollen, we will slightly modify this formula. In order to mimic Definition \ref{DefQualityNectar}, we need a reference flower that only gives pollen to bees. For any other flower $f$ that produces pollen, the quality of its pollen relative to the reference flower is
\begin{align*}
\frac{\tilde{q}_f}{\tilde{q}_0},
\end{align*}
which represents the ratio of the average lifespan of bees fed with pollen from flower $f$ relative to the average lifespan of bees fed with the reference flower \cite{PollenNutrition}. The tilde on $\tilde{q}$ means that we are considering pollen.
\newline
Now the quality of the pollen is
\begin{equation}
\boxed{\tilde{\mathcal{Q}}_f\left(\tau\right)=\overbrace{\tau\frac{\tilde{q}_f}{\tilde{q}_0}}^{\tilde{\mathcal{Q}}_{f,costless}\left( \tau \right)} -\overbrace{\xi\times \left(\frac{d}{d_{max}} +\sqrt[n]{\frac{\rho_{crit,f}}{\rho_f}}\right)}^{\tilde{\mathcal{Q}}_{f,cost}},}
\end{equation}
where $\xi$ represents the change in the aerodynamics of the bee due to the pollen it is holding, and can be replaced by $1$ as a first approximation. Notice that in the previous quality definition we \underline{did not} substract a term corresponding to the honey the bee uses to aggregate the pollen balls it carries. Indeed, this honey is not lost by the colony as it is mixed to pollen that will feed other bees.
\newline
To define $\tilde{\eta}_f \left( \tau \right)$ , $\tilde{\eta}_{f,costless} \left( \tau \right)$ and $\tilde{\eta}_{f,cost}$  we only need to replace $\mathcal{Q}_f$ by $\tilde{\mathcal{Q}}_f \left( \tau \right)$, $\tilde{\mathcal{Q}}_{f,costless} \left( \tau \right)$ or $\tilde{\mathcal{Q}}_{f,cost}$ in Definition \ref{DefEfficient}.
\newline
Now let us explain how $\tau$ is defined.
\begin{algo} \label{ComputeTau} A)\underline{In the case where $\mathcal{R}_{hive}<\mathcal{R}_{target}$,}
\newline
We establish a list of nectar resources ordered by efficiency $\eta_f$ and we express the same list for pollen by using efficiency $\tilde{\eta}_{f,costless}\left(\tau\right)$ to fix the order on the resources (that does not depend on the choice of $\tau$).
\newline
Since $\mathcal{R}_{hive}<\mathcal{R}_{target}$, there is not enough honey in the hive so :
\newline $\bullet$ we appoint enough bees on the pollen resources where the bees are most efficient in order for the bees to bring back a total income corresponding to a fixed minimal value (we use $\tilde{\eta}_{f,costless}\left(\tau\right)$ instead of $\tilde{\eta}_{f}\left(\tau\right)$ in the case $\mathcal{R}_{hive}<\mathcal{R}_{target}$. Indeed the bees must necessarily fetch the minimal pollen income so they may apply the "whatever it takes" policy to do so),
\newline $\bullet$ all the other forager bees are sent to the nectar resources where they are the more efficient to maximize nectar income.
\newline
The choice of the nectar resources to which the nectar foragers are sent corresponds to a value $\mathcal{\eta}_{cut}$ above which the nectar resources are considered useful. $\tau$ is then defined to satisfy
\begin{equation*}
\boxed{\tilde{\eta}_{cut,costless}(\tau) =\eta_{cut} +\frac{\sum\limits_{\substack{\text{useful}\\\text{pollinic }f}} \phi_f\tilde{\eta}_{f,cost}}{\sum\limits_{\substack{\text{useful}\\\text{pollinic }f}}\phi_f}.}\\
\end{equation*}
In the previous formula, we consider that bees that forage on different pollinic resources "compensate" each other's cost in order for them all to pay an equal cost, which makes sense because we are not focusing on the scale of individual bees but on the scale of the hive. We moreover recall that the number of bees $\phi_f$ that forage on flower $f$ is given by formula \ref{NumberBeesresource}.
\begin{center}
\includegraphics[scale=.5]{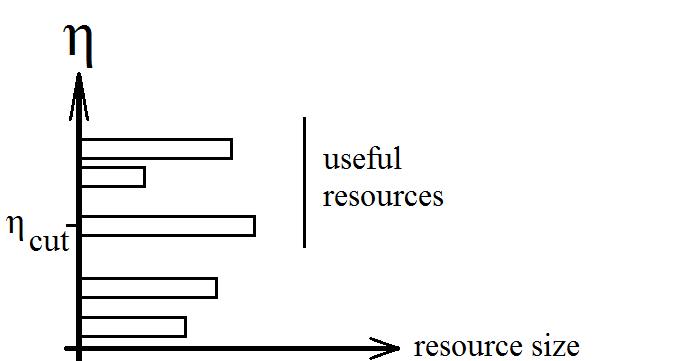}
\end{center}
B)\underline{In the case where $\mathcal{R}_{hive}>\mathcal{R}_{target}$,}
\newline
there is not enough pollen in the hive in comparison to honey.
\newline
As was done in the previous case, we start by attributing the task of collecting a minimal amount of the non-scarce resource (here nectar) to a small number of bees. This minimal amount compensates energy losses  infuenced by such parameters as the outside temperature and the number of bees in the hive.
\newline $\bullet$ We appoint enough bees on the best resources (according to the $\eta_{f}$ metric) to bring back the right amount of nectar. This defines
\begin{align*}
\eta_{cut}=\min \left\lbrace \eta_f \vert f \text{ is visited by bees} \right\rbrace .
\end{align*} 
\newline The foragers that remain after some bees were sent to nectar resources must be sent to bring pollen back. For any choice of $t\in \mathbb{R}^{+*}$, by considering the efficiency metric $\tilde{\eta}_{f}(t)$, we rate the floral resources and send the bees to the best ones.  $\tilde{\eta}_{cut}(t)$ can be defined the same way $\eta_{cut}$ was defined. As a function $\tilde{\eta}_{cut}$ is equal to zero on $\left[ 0,t_0 \right]$ and is continuous, piecewise affine and strictly increasing above $t_0$. Moreover, if $t_0<+\infty$ there exists a value $t_1$ above which $\tilde{\eta}_{cut}$ is affine and strictly increasing.
\newline
We will justify all the properties satisfied by $\tilde{\eta}_{cut}$ by building its graph by using a "colorization" of the plan after we finish the description of the algorithm to build $\tau$.
\newline
i)\underline{In the case where $t_0<+\infty$} there exists $\tau_1$ such $\tilde{\eta}_{cut}\left( \tau_1 \right)=\eta_{cut}$. $\tau_1$ corresponds via the $\tilde{\eta}_{f}$ functions to a certain way to distribute bees on resources and the cost of fetching pollen can now be estimated. This estimation of the cost can now be taken into account into the minimal quantity of nectar that must be brought back to the hive and we can start the algorithm (case B) ) over again with this new data. This gives us a new potential exchange rate $\tau_2$ but since more bees went to get nectar, less were going to pollen resources and the estimation of the cost of fetching pollen is once again corrected in the other direction. We go on and compute $\tau_3$, $\tau_4$, ... Notice that by construction $\left( \tau_{2i} \right)$ and $\left( \tau_{2i+1} \right)$ are adjacent. The exchange rate is then defined as
\begin{align*}
\tau =\lim\limits_{i\rightarrow \infty} \frac{\tau_{i}+\tau_{i+1}}{2}.
\end{align*}
\newline
ii)\underline{The case where $t_0=+\infty$} corresponds to scarce resources where some forager bees may risk unemployement if the hive goes at full speed. A strategy the beehive may adopt is to adapt the previous algorithm by substituting quality $\mathcal{Q}$ to efficiency $\eta$.
\end{algo}
We still need to explain how to easily construct the graph of $\tilde{\eta}_{cut}$ in the previous algorithm (case B) by coloring the plan.
\newline
For every floral resource $f$, we draw the graph of $\tilde{\eta}_f$ which is an affine function of $\tau$. We color its hypograph with a level of grey proportional to the number of bees $\phi_f$ that it takes to  exhaust it and that is given by formula \ref{NumberBeesresource}. Where there is a superposition of levels of grey, we add these levels of grey. The graph of the function $\tilde{\eta}_{cut}$ can then be seen as the boundary between a level of grey that represents less than the number of bees that are affected to pollen resources and a level of grey that represents more than this number of bees.
\begin{center}
\includegraphics[scale=.5]{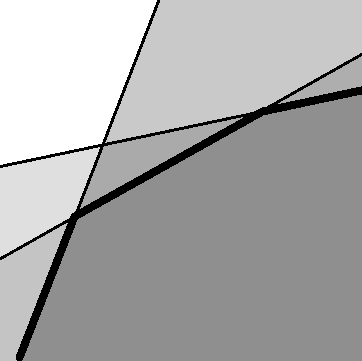}
\end{center}

\section{Predators}
\subsection{Predation during foraging activities}
In practice, in a certain region, the maximal distance at which bees may forage for nectar is not the theoretical $d_{max}$, even if the floral resources are scarce and the flower that is at $d_{max}$ is the reference flower with optimal $q_f$. Indeed, $\frac{q_0}{d_{max}}$ represented the energetical cost of flying at a certain distance but there is an additionnal cost that is not included in $\frac{q_0}{d_{max}}$ that is the cost of predation on the bees that go from the hive to a floral resource. To include this cost that depends on the local fauna, we modify $d_{max}$ and transform it into $d_{max, local}$ so that the total cost of flying to a floral resource per distance unit is
\begin{align*}
\frac{q_0}{d_{max,local}},
\end{align*}
Which means that the cost of predation per distance unit can be expressed as
\begin{align*}
q_0\left(\frac{1}{d_{max,local}}-\frac{1}{d_{max}}\right).
\end{align*}
When one forager bee dies by predation, the energy-equivalent cost for the colony is
\begin{align*}
\alpha+\tau\frac{\tilde{\alpha}q_0}{\tilde{q}_0}\times\overbrace{\frac{\mathcal{L}_{forager}}{2 \mathcal{L}_{average}}}^{\substack{\text{pro rata}\\ \text{remaining}\\ \text{lifespan}}},
\end{align*}
where
\newline $\bullet$ $\alpha$ is the energy used to "build" one bee which appears as a coefficient in Equation \ref{bilanenergie},
\newline $\bullet$ $\tilde{\alpha}$ is the quantity of pollen used to "build" one bee (see \cite{Pollen1}),
\newline $\bullet$ $\mathcal{L}_{forager}$ is the average time a bee stays a forager and $\mathcal{L}_{average}$ is the average total lifespan of one bee. Data on these lifespans can be found for instance in \cite{RFID2}.
\newline
The rate of death by predation per time unit during flights to floral resources can be estimated to be 
\begin{equation}
\boxed{\mathcal{P}_{flight}= V \times \frac{ q_0\left( \frac{1}{d_{max,local}}-\frac{1}{d_{max}} \right) }{ \alpha \left(  1+\tau\frac{\tilde{\alpha}q_0}{\alpha\tilde{q}_0}\times\frac{\mathcal{L}_{forager}}{2 \mathcal{L}_{average}} \right)  } . }
\end{equation}
Similarily, we can write the rate of predation per time unit on the foraging site $f$ as 
\begin{equation}
\boxed{\mathcal{P}_{foraging}= \frac{1}{\frac{1}{v_{hop}}+\frac{\beta_f \sqrt[n]{\rho_f}}{k_n}} \times \frac{ q_0\left( {\sqrt[n]{\rho_{crit,local}}-\sqrt[n]{\rho_{crit}}} \right) }{ \alpha k_n \left(  1+\tau\frac{\tilde{\alpha}q_0}{\alpha\tilde{q}_0}\times\frac{\mathcal{L}_{forager}}{2 \mathcal{L}_{average}} \right)  } . }
\end{equation}

\subsection{Hornets}
A bee is to a hornet what pollen is to a bee. In this analogy, the beekeeper can be considered as a predator from the perspective of the hornet and if this predation is too costly for him, he will not venture near the beehive.
\section{Local bees}
Throughout this text, there appear many clear reasons why bees that are adapted to a local environnement should be preserved :
\newline $\bullet$ To optimize the heating of the hive, a coefficient $\nu$ is chosen that implies a temperature inside the bee cluster given by formula \ref{TemperatureCluster}. The colder the winter is supposed to be, the higher $\nu$ must be, which costs more energy. Local bees are able to find the correct value for $\nu$.
\newline $\bullet$ Also in order to survive during winter, local bees know how much honey to stock, which means that they compute $\mathcal{R}_{target}$ correctly.
\newline $\bullet$ Local bees know how to evaluate $q_f$ and $\lambda_f$ for every flower spieces around them.
\newline $\bullet$ The seasonal adjustments of the hive between $\mathcal{R}_{hive}>\mathcal{R}_{target}$ and $\mathcal{R}_{hive}<\mathcal{R}_{target}$ have to be in phase with local flora. Local bees have found seasonal cycles that harmoniously fit in the flora.
\newline $\bullet$ Local bees know local predators.

\vskip.5cm

\end{document}